\documentclass[11pt]{article}
\usepackage{amsthm,amsmath,latexsym,amssymb}

\newcommand{\bP}{{\rm |\kern-.15em P}}
\newcommand{\Q}{\kern.3em\rule{.07em}{.65em}\kern-.3em{\rm Q}}
\newcommand{\R}{{\rm I\kern-.15em R}}
\newcommand{\D}{{\rm |\kern-.15em D}}
\newcommand{\h}{{\rm |\kern-.15em H}}
\newcommand{\C}{\kern.3em\rule{.07em}{.65em}\kern-.3em{\rm C}}
\newcommand{\T}{{\rm T\kern-.35em T}}

\theoremstyle{plain}
\newtheorem{theorem}{Theorem}[section]

\newtheorem{proposition}[theorem]{Proposition}

\theoremstyle{definition}
\newtheorem{definition}[theorem]{Definition}
\newtheorem{example}[theorem]{Example}

\theoremstyle{remark}
\newtheorem{remark}[theorem]{Remark}

\begin{document}
\title{On the structure of the semigroup of entire \'etale mappings}
\author{Ronen Peretz}
 
\maketitle

\section{Introduction}
This is a paper in analysis. In fact it deals with entire functions in a single complex variable which
are local homeomorphisms of the complex plane. Also a central tool used is the so called composition operators.
However, the motivation for the proposed research originates in the two-dimensional Jacobian Conjecture.
The Jacobian Conjecture is one of the most famous and prominent conjectures in algebraic
geometry \cite{bcw,e,s}. It is mentioned in the list of open problems prepared by Steve Smale \cite{s} that intends to orient
mathematical research in this century (parallel to the famous list of problems prepared by David
Hilbert at the beginning of the previous century). Thus, any significant progress in this problem is expected
to make an important contribution to mathematical science. If it is true, then the Jacobian Conjecture gives a remarkably simple necessary and sufficient
condition on a polynomial mapping $F:\,\mathbb{C}^n\rightarrow\mathbb{C}^n$ to be an automorphism \cite{e1,e}.
An automorphism is an invertible mapping $\mathbb{C}^n\rightarrow\mathbb{C}^n$ (i.e. an injective and a surjective mapping whose inverse 
$F^{-1}:\,\mathbb{C}^n\rightarrow\mathbb{C}^n$ is also polynomial). These polynomial mappings are precisely the morphisms 
that preserve the  algebro-geometric affine structure of $\mathbb{C}^n$, which is a central theme in affine algebraic
geometry. This is one instance of the importance of the conjecture.

Let us denote by $J(F)$ or by $J_F$ the Jacobian matrix of the mapping $F$. If $F$ is invertible,
then in particular it is locally invertible and hence by the Inverse Mapping Theorem the determinant
of its Jacobian matrix, $\det J(F)$, does not vanish at any point of $\mathbb{C}^n$. However,
the fact that $F$ is a polynomial mapping implies that $\det J(F)$ is a polynomial over the complex
field $\mathbb{C}$ and hence the Fundamental Theorem of Algebra implies that $\det J(F)$ must be
a non-zero constant. This elementary argument proves that a necessary condition for a polynomial
mapping $F$ to be an automorphism of $\mathbb{C}^n$, i.e. $F\in {\rm Aut}(\mathbb{C}^n)$, is
that $\det J(F)\in\mathbb{C}^{\times}$. The Jacobian Conjecture speculates the validity of the 
inverse statement. Thus \underline{{\bf the Jacobian Conjecture}} is
$$
F\in {\rm Aut}(\mathbb{C}^n)\Leftrightarrow \det J(F)\in\mathbb{C}^{\times}.
$$
This is true for dimension $n=1$ but it is wide open for dimension $n\ge 2$. 
The original version of the conjecture was stated by Ott Keller in 1939 \cite{k}.
Keller worked over the integers $\mathbb{Z}$ and with polynomial mappings $F:\,\mathbb{Z}^n
\rightarrow\mathbb{Z}^n$ that satisfy the unimodular Jacobian condition $\det J(F)=1\,\,{\rm or}\,\,(-1)$.
Since that time, much research has been done on the so-called Keller Problem, giving rise to great many
beautiful ideas and theories. Notable results are the degree reduction theorems \cite{bcw,d1,d2,y}. These are based
on $K$-theoretic principles that allow us, for example, to reduce the proof of the conjecture to the seemingly
simple case of mappings of the form $F=(X_1,\ldots,X_n)+H(X_1,\ldots,X_n)$, where $(X_1,\ldots,X_n)$
is the identity mapping and $H(X_1,\ldots,X_n)$ is a cubic homogeneous mapping. However, one should
prove the Jacobian Conjecture for such cubic mappings in an arbitrary dimension. We also note that the conjecture is known to be true
in any dimension and for mappings of degree $\deg F\le 2$, \cite{bcw,e}. So it seems that "we are almost there".
Yet the degree $3$ case seems to be out of our reach at least for the present.

Some experts in this field tend to 
believe that the general conjecture ($n\ge 2$) might be false but that the two-dimensional case is
possibly true. We mention here two pivotal two dimensional results. The first is a theorem of Moh \cite{moh}, 
which asserts the validity of the conjecture for $n=2$ and degree $d=\deg F\le 100$ or so. This 
is a difficult result. The second result is the ingenious counterexample of Pinchuk \cite{pin} to the so-called
Real Jacobian Conjecture \cite{r}; namely, if $F:\,\mathbb{R}^2\rightarrow\mathbb{R}^2$ is a real polynomial
mapping that satisfies the real Jacobian condition $\det J_F(X,Y)\ne 0\,\,\forall\,(X,Y)\in\mathbb{R}^2$,
then $F^{-1}$ exists. We note that in this case $\det J_F$ need not be a non-zero constant, but that the
conclusion is also weaker, namely, $F^{-1}$ does not need to be a polynomial mapping. This natural real version
of the original Jacobian Conjecture was open for sometime until in 1993 Pinchuk found a counterexample.
His original clever construction gave rise to a degree $35$ counterexample, but almost immediately it
was reduced to a degree $25$ counterexample. It is interesting to mention that the minimal degree of
a counterexample is still not known, but we do have some lower bounds ($d\ge 7$).

Another approach to solving the Jacobian Conjecture is via a thorough analysis of the structure of the semigroup
of all the normalized \'etale mappings on $\mathbb{C}^n$.This approach is outlined in the papers of
Kambyashi \cite{k1,k2,k3} and of the author \cite{rp0,rp1}.
The set of all the normalized \'etale mappings on $\mathbb{C}^n$ of degree $d$ or less can naturally be
parametrized on a finite dimensional algebraic space,  the so called Jacobian variety, $J(n,d)$, of degree
$d$. Taking the union of these varieties gives us an object $J(n)$, which is known in the literature as
an ind-variety. There is a canonical way to define a topology on $J(n)$. We call $J(n)$ the Jacobian variety
(of $\mathbb{C}^n$). As shown in the above papers, this structure might lead us to make a progress on the 
Jacobian Conjecture. Moreover, it is apparent (see \cite{rp}) that using this theory to solve the Jacobian
Conjecture is intimately related to the singular loci $S(J(n,d))$ of the finite dimensional varieties. 
Equally relevant to the conjecture might be the singular locus $S(J(n))$ of the
total ind-variety. However, it is still a challenge in the theory of ind-varieties to give a satisfactory
definition of the smooth points. An early attempt by Shafarevich \cite{sh1,sh2} did not work well,
as indicated by the clever example of Totatro \cite{k1,k2}. Still it is possible that Shafarevich's
definition will work well for our particular ind-variety $J(n)$. A good source to consult for singular
ind-varieties is the book \cite{kumar} of Kumar. 

In our proposed research we hope to overcome that
difficulty of attacking the Jacobian Conjecture via the structure of ind-varieties, at least in the 
two-dimensional case. We plan to do that by imposing a different structure on the family of the two-dimensional
\'etale mappings. This structure is a geometric structure, as opposed to the algebraic ind-variety structure.
Our proposed structure is a fractal structure where the role of the degree filter in the ind-variety structure,
will be taken by self-similarity structure that is induced by the binary operation of the semigroup, i.e. the composition
operators. We have two of these operators, a left composition operator and a right composition operator. They
have different properties!

We hope to be able (in the future) to carefully analyze this structure and to tie the Jacobian Conjecture
in dimension two to certain Zeta functions, thereby invoking a powerful arithmetic machinery to handle
the two dimensional Jacobian Conjecture. Let us denote by ${\rm et}(\mathbb{C}^2)$ the semigroup of two dimensional
Keller mappings. We would like to prove something like the following: \\
\\
a) That there exists an infinite index set $I$, and a family of mappings indexed by $I$, $\{ F_i\,|\,i\in I\}
\subset {\rm et}(\mathbb{C}^2)$ such that 
$$
{\rm et}(\mathbb{C}^2)={\rm Aut}(\mathbb{C}^2)\cup\bigcup_{i\in I} R_{F_i}({\rm et}(\mathbb{C}^2)),
$$
where if $i\ne j$ then $R_{F_i}({\rm et}(\mathbb{C}^2))\cap R_{F_j}({\rm et}(\mathbb{C}^2))=\emptyset$. \\
b) That the parallel representation to the representation described in (a) above holds true, 
this time with respect to the left composition operators $L_{G_j}$. 

These two claims will be the basis for a fractal structure on ${\rm et}(\mathbb{C}^2)$ because the pieces $R_{F_i}({\rm et}(\mathbb{C}^2))$
are similar to each other in the sense that they are homeomorphic, and we further have the scaling property of
self-similarity, namely $R_F({\rm et}(\mathbb{C}^2))$ is homeomorphic to its proper subspace \\
$R_{G\circ F}({\rm et}(\mathbb{C}^2))$ that is homeomorphic to its proper subspace $R_{H\circ G\circ F}({\rm et}(\mathbb{C}^2))$
etc . This is the right place to remark that the purpose of the current paper is to start and develop the parallel theory for entire functions in one complex
variable. Results in this setting will hint that there are hopes to accomplish the above objective. 

But our construction goes even further than this! The structure we propose to investigate is far more general than the structure
of ind-varieties outlined  in Kambyashi's papers \cite{k1,k2,k3}. The filtration in the ind-variety 
structure is an algebraic degree filtration. Thus it makes sense only in the algebraic context of polynomial
mappings (as in the Jacobian Conjecture). However, the geometric structure that we suggest relies only
on the fact that the underlined semigroup is a semigroup of local diffeomorphisms and is not necessarily algebraic.
Thus, our theory applies to a much broader context of mappings. In addition, in that universe, we hope to be able
to relate the semigroup structure that we suggest to Hausdorff measures and dimensions and to arithmetic. Thus
we expect that our theory will also make a contribution in the opposite direction; namely, starting from
a relatively well-understood semigroup of local diffeomorphisms, we will set out to conclude non-trivial theorems
in the theories of Hausdorff measures, in fractals and in arithmetic.

The purpose of this paper is to deal with a semigroup very different from the semigroup of \'etale
polynomial mappings on $\mathbb{C}^2$; namely, we will start to investigate the normalized semigroup
of entire functions in one complex variable that has a nonvanishing first order derivative. This semigroup
is better understood than the algebraic one. For example, we can (and will) invoke Nevanlinna's value distribution
theory for entire functions in one complex variable. In the future we will try to understand what features are common to
this analytic semigroup and to the algebraic semigroup that is suitable for the two-dimensional
Jacobian Conjecture. A few surprises are to be expected. The reason is that our intuition will be 
built upon the analytic semigroup and so will lead us to expect similar properties in the algebraic setting.
However, we might find out that not always will the algebraic semigroup have the same properties as the 
analytic semigroup.

\section{The semigroup of the entire local homeomorphisms in one complex variable, 
normalized at the origin}
\begin{definition}
$$
{\rm elh}(\mathbb{C})=\{ f\,:\,\mathbb{C}\rightarrow\mathbb{C}\,|\,f\,{\rm is}\,{\rm entire},\,
\forall\,z\in\mathbb{C}\,\,f'(z)\ne 0,\,\,f'(0)=1\},
$$
$$
{\rm Aut}(\mathbb{C})=\{ z+a\,|\,a\in\mathbb{C}\}.
$$
\end{definition}

\begin{proposition}
{\rm 1)} $({\rm elh}(\mathbb{C}),\circ)$ is a semigroup with an identity ($\circ $ is composition). \\
{\rm 2)} If $\{f_n\}$ is a sequence in ${\rm elh}(\mathbb{C})$ that converges to $f$ uniformly on compact
subsets of $\mathbb{C}$, then $f\in {\rm elh}(\mathbb{C})$. \\
{\rm 3)} If $\{f_n\}$ is a sequence in ${\rm elh}(\mathbb{C})$ that satisfies the Cauchy condition
uniformly on compact subsets of $\mathbb{C}$, then $\lim f_n\in {\rm elh}(\mathbb{C})$. \\
{\rm 4)} $({\rm Aut}(\mathbb{C}),\circ)$ is a closed group in $({\rm elh}(\mathbb{C}),\circ)$ with respect
to the topology of local uniform convergence (i.e. convergence uniformly on compact subsets of $\mathbb{C}$). \\
{\rm 5)} $f\in {\rm elh}(\mathbb{C})-{\rm Aut}(\mathbb{C})\,\Leftrightarrow A(f)\ne\emptyset$ (See remark 2.3 below). \\
{\rm 6)} $\forall\,f\in {\rm elh}(\mathbb{C}),\,\forall\,z\in\mathbb{C}$ the fiber $f^{-1}(z)$ is a discrete
subset of $\mathbb{C}$. \\
{\rm 7)} $\forall\,f\in {\rm elh}(\mathbb{C}),\,|\mathbb{C}-f(\mathbb{C})|<2$.
\end{proposition}

\begin{remark}
We recall that the asymptotic variety of $f$, $A(f)$, mentioned in part 5 of the proposition is defined as 
follows:
$$
A(f)=\{\alpha\in\mathbb{C}\,|\,\exists\,\gamma:\,[0,\infty)\rightarrow\mathbb{C}\,\,{\rm continuous}\,\,
{\rm such}\,\,{\rm that}\,\,\lim_{t\rightarrow\infty}|\gamma(t)|=\infty
$$
$$
{\rm and}\,\,\lim_{t\rightarrow\infty} f(\gamma(t))=\alpha\}.
$$
We don't use the sequential asymptotic variety $A_S(f)$:
$$
A_S(f)=\{\alpha\in\mathbb{C}\,|\,\exists\,\{ z_n\}\subseteq\mathbb{C},\,\,{\rm such}\,\,{\rm that}\,\,\lim |z_n|=\infty
\,\,{\rm and}\,\,\lim f(z_n)=\alpha\}.
$$
The reason will soon be explained. We clearly have: $A(f)\subseteq A_S(f)$. Also $A(f)$ contains the set
of all the $f$-Picard values, $\mathbb{C}-f(\mathbb{C})$.
\end{remark}
\noindent
{\bf A proof of Proposition 2.2.} \\
1) By the chain rule, $f,g\in {\rm elh}(\mathbb{C})\Rightarrow f\circ g\in {\rm elh}(\mathbb{C})$. The identity
is ${\rm id}(z):\equiv z$. \\
2) By Cauchy's Theorem, $f=\lim f_n$ is an entire function. By Cauchy's estimate, $f'=\lim f'_n$ locally
uniformly. By Hurwitz Theorem either $f'(z)\equiv 0$ or $f'$ never vanishes. By the normalization at
the origin, $f'(0)=1$ and hence $f\in {\rm elh}(\mathbb{C})$. \\
3) The assumption is that $\forall\,K\subset\mathbb{C}$ a compact, and $\forall\,\epsilon>0$ $\exists\,
N_{\epsilon}(K)$ such that $\forall\,n,m\in\mathbb{Z}^+$, we have
$$
n,m> N_{\epsilon}(K)\Rightarrow \forall\,z\in K,\,\,|f_n(z)-f_m(z)|<\epsilon.
$$
So $\{f_n\}$ converges locally uniformly on $\mathbb{C}$ and the conclusion follows by
part 2. \\
4) Clearly $({\rm Aut}(\mathbb{C}),\circ )$ is a group. If $\{z+a_n\}$ converges locally uniformly
then $a=\lim a_n\in\mathbb{C}$ exists and $\lim(z+a_n)=z+a\in{\rm Aut}(\mathbb{C})$. \\
5) By Hadamard's Theorem, \cite{h}, $f\in {\rm elh}(\mathbb{C})-{\rm Aut}(\mathbb{C})\Leftrightarrow A(f)\ne
\emptyset$. \\
6) If $z_0\in\mathbb{C}$, $f\in{\rm elh}(\mathbb{C})$ and $f^{-1}(z_0)$ is not a discrete subset of $\mathbb{C}$,
then $\exists\,\{w_n\}\subseteq f^{-1}(z_0)$ of different points, such that $w_0=\lim w_n$ exists. Hence
$\forall\,n\in\mathbb{Z}^+$, $f(w_n)=z_0$ and by the permanence principle $f(z)\equiv z_0$ which
contradicts $f\in{\rm elh}(\mathbb{C})$. \\
7) This is the Picard's Little Theorem. $\qed $ \\
\\
We now explain why we do not use $A_S(f)$.

\begin{proposition}
If $f(z)$ is an entire non-polynomial function, then $A_S(f)=\mathbb{C}$.
\end{proposition}
\noindent
{\bf Proof.} \\
Let $z_0\in\mathbb{C}$ and let $\{ w_n\}\subseteq f(\mathbb{C})$ satisfy $\lim w_n=z_0$. Since $f$ is
non-polynomial it has infinitely many regular maximal domains $\Omega_n$ that tile up the complex plane,
\cite{rp}. We choose $\forall\,n\in\mathbb{Z}^+$, $z_n\in\Omega_n$ such that $w_n=f(z_n)$. By passing to
a subsequence, if necessary, we may assume that $\lim |z_n|=\infty$. Hence $z_0\in A_S(f)$. $\qed$ \\

\section{The right and the left mappings on ${\rm elh}(\mathbb{C})$}

\begin{definition}
Let $f\in {\rm elh}(\mathbb{C})$. The right mapping on ${\rm elh}(\mathbb{C})$, induced by $f$ is defined by:
$$
\begin{array}{l} R_f:\,{\rm elh}(\mathbb{C})\rightarrow {\rm elh}(\mathbb{C}) \\ R_f(g)=g\circ f.\end{array}
$$
The left mapping on ${\rm elh}(\mathbb{C})$, induced by $f$ is defined by:
$$
\begin{array}{l} L_f:\,{\rm elh}(\mathbb{C})\rightarrow {\rm elh}(\mathbb{C}) \\ L_f(g)=f\circ g.\end{array}
$$
\end{definition}

\begin{definition}
Let $g\in {\rm elh}(\mathbb{C})$. We will denote by $T_f(g)$ the set of all the finite asymptotic tracts of $g(z)$. Thus
formally:
$$
T_f(g)=\{\gamma:\,[0,\infty)\rightarrow\mathbb{C}\,|\,\gamma\,{\rm continuous},\,\lim_{t\rightarrow\infty}|\gamma(t)|
=\infty\,{\rm and}\,\lim_{t\rightarrow\infty} g(\gamma(t))\in\mathbb{C}\}.
$$
\end{definition}
\noindent
We did not include in $T_f(g)$ the asymptotic tracts of $g$ of the infinite asymptotic value, i.e. those tracts
$\gamma:\,[0,\infty)\rightarrow\mathbb{C}$ which are continuous and satisfy $\lim_{t\rightarrow\infty}|\gamma(t)|=\infty$
and $\lim_{t\rightarrow\infty} |g(\gamma(t))|=\infty$. We collect together asymptotic tracts of $g$ according to
the following criterion.

\begin{definition}
Two finite asymptotic tracts $\gamma,\beta\in T_f(g)$ of $g\in {\rm elh}(\mathbb{C})$ are said to be 
equivalent if: \\
a) $\lim_{t\rightarrow\infty} g(\gamma(t))=\lim_{t\rightarrow\infty} g(\beta(t))$. \\
b) The finite tracts $\gamma$ and $\beta$ are homotopic via a homotopy that fixes $\infty$. Thus 
$\exists\,H(s,t):\,[0,1]\times [0,\infty)\rightarrow\mathbb{C}$ continuous in $(s,t)$ and satisfying
$H(0,t)=\gamma(t)$, $H(1,t)=\beta(t)\,\,\forall\,t\ge 0$ and 
$$
\left\{\begin{array}{l} \lim_{t\rightarrow\infty}|H(s,t)|=\infty\,\forall\,0\le s\le 1, \\
\lim_{t\rightarrow\infty}f(H(s,t))=\lim_{t\rightarrow\infty}g(\gamma(t)).\end{array}\right.
$$ 
It is standard to conclude that the above relation on $T_f(g)$ is an equivalence relation. We denote
by $T_{0f}(g)$ the set of all the equivalence classes.
\end{definition}
\noindent
\noindent
It is natural to ask if there is a canonical way to choose representatives of the classes in $T_{0f}(g)$.
We would be happy to obtain such a representative by a construction similar to the resolution process
that was used in $\mathbb{C}[X,Y]^2$ (See \cite{p3}). However, in trying to get the parallel of the asymptotic identities
we find that the situation in ${\rm elh}(\mathbb{C})$ is harder. To explain we now prove the following proposition.

\begin{proposition}
Let $f(z)$ be a non-constant entire function. \\
1) If $g(z)\in H(0<|z-z_0|\le r)-H(|z-z_0|\le r)$, then $f(g(z))\in H(0<|z-z_0|\le r)-H(|z-z_0|\le r)$. \\
2) If $h(z,w)$ has a non-removable singularity at $(z_0,w_0)$ then $f(h(z,w))$ can not be holomorphic in
a neighborhood of $(z_0,w_0)$.
\end{proposition}
\noindent
The symbol $H(\Omega)$ for an open non-empty $\Omega\subseteq\mathbb{C}$ stands for the family of holomorphic
functions defined on $\Omega$. \\
\\
{\bf Proof.} \\
1) By the assumption $z_0$ is an isolated non-removable singularity of $g(z)$. Hence either $z_0$ is an essential
singularity of $g(z)$ or it is a pole of $g(z)$. \\
{\bf Case 1:} $z_0$ is an essential singularity of $g(z)$. Then by the Picard Theorem, $g(\{0<|z-z_0|\le r\})=\mathbb{C}$
or it is $\mathbb{C}-\{a\}$ for a fixed $a\in\mathbb{C}$. If $f(g(z))\in H(|z-z_0|\le r)$ then $f(g(\{0<|z-z_0|\le r\})$
had a compact closure and hence $f(\mathbb{C})$ or $f(\mathbb{C}-\{a\})$ had a compact closure. By the Liouville Theorem
$f(z)$ had to be a constant, which contradicts our assumption on $f(z)$. \\
{\bf Case 2:} $z_0$ is a pole of $g(z)$. In that case $\lim_{z\rightarrow z_0}g(z)=\infty$ and $g(\{0<|z-z_0|\le r\})$ contains
the complement of some disc, say $\{R<|w|\}$. If $f(g(z))\in H(|z-z_0|\le r)$ then $f(\{R<|w|\})$ had a compact
closure and so also had $f(\mathbb{C})$. By Liuoville Theorem $f(z)$ had to be a constant, which is a contradiction. \\
2) In this case we can either fix $w=w_0'$ and consider $f(h(z,w_0'))$ or fix $z=z_0'$ and consider $f(h(z_0',w))$ and
then use part 1. $\qed$ \\
\\
So we can not expect naive asymptotic identities in the holomorphic one variable case, which will resemble the
two variable polynomial case. For general entire functions the asymptotic variety $A(f)$ can be pretty wild. For 
instance it could be any finite subset of $\mathbb{C}$, and it could be as large as the whole space, $\mathbb{C}$.

\begin{proposition}
If $f,g\in {\rm elh}(\mathbb{C})$ then $T_f(g)\subseteq T_f(f\circ g)$, $f(A(g))\subseteq  A(f\circ g)$
\end{proposition}
\noindent
{\bf Proof.} \\
$\gamma\in T_f(g)\Rightarrow \gamma:\,[0,\infty)\rightarrow\mathbb{C}$ is continuous, $\lim_{t\rightarrow\infty}|\gamma(t)|=\infty$
and $\lim_{t\rightarrow\infty}g(\gamma(t))\in\mathbb{C}$ $\Rightarrow\lim_{t\rightarrow\infty}f(g(\gamma(t)))\in\mathbb{C}
\Rightarrow\lim_{t\rightarrow\infty}(f\circ g)(\gamma(t))\in\mathbb{C}\Rightarrow\gamma\in T_f(f\circ g)$. \\
$a\in f(A(g))\Rightarrow\exists\,\gamma\in T_f(g)$ such that $a=\lim_{t\rightarrow\infty}f(g(\gamma(t)))\Rightarrow
\exists\,\gamma\in T_f(f\circ g)$ such that $a=\lim_{t\rightarrow\infty} (f\circ g)(\gamma(t))\Rightarrow
a\in A(f\circ g)$. $\qed $ \\
\\
This proposition tells us that left composition of functions in ${\rm elh}(\mathbb{C})$ does nor decrease $T_f(g)$ ($g$ the right function), 
the set of all the finite asymptotic tracts of the right factor and consequently it does not decrease the left image of 
its finite asymptotic values. We ask, when in Proposition 3.5 we have a strict containment? This happens
exactly when $\exists\,\gamma\in T_f(f\circ g)-T_(g)$. This implies that $a=\lim_{t\rightarrow\infty}(f\circ g)(\gamma(t))\in\mathbb{C}$
but $\lim_{t\rightarrow\infty}g(\gamma(t))$ does not exist. If by "does not exist" we insist that 
$\lim_{t\rightarrow\infty}g(\gamma(t))$ does not exist as a finite number, then of course we can take $g(z)=z$
and $\gamma\in T_f(f)$. However, if we make no distinction between finite limits and an infinite limit
then this example no longer works. Thus the following is worth mentioning.

\begin{proposition}
If $f,g\in {\rm elh}(\mathbb{C})$, then $\gamma\in T_f(f\circ g)$ implies that either $\gamma\in T_f(g)$ or else
$\lim_{t\rightarrow\infty} g(\gamma(t))=\infty$.
\end{proposition}
\noindent
{\bf Proof.} \\
Of course, it is enough to take $\gamma\in T_f(f\circ g)-T_f(g)$ and prove that $\lim_{t\rightarrow\infty}g(\gamma(t))=\infty$.
If this is not the case then $\lim_{t\rightarrow\infty}g(\gamma(t))$ does not exist in the broad sense (i.e. not
as a finite number nor $\infty$). Hence we have two sequences $t_n,s_n\in [0,\infty)$ so that $\lim t_n=\lim s_n=\infty$ and
$a=\lim g(\gamma(t_n))$, $b=\lim g(\gamma(s_n))$ exist and are not equal to one another. But by our assumption
$\gamma\in T_f(f\circ g)$ and so $f(a)=f(b)$. We claim that there exists a sequence $\{ a_n\}\subseteq\mathbb{C}$ of
distinct points for which $\lim a_n=a$ and $\forall\,n,\,f(a_n)=f(a)$. If indeed this is true then by the permanence 
principle $f(z)\equiv f(a)$ a constant. This contradicts the assumption that $f\in {\rm elh}(\mathbb{C})$. To prove the
existence of the sequence $\{ a_n\}$ we prove the following: $\forall\,0<r<|b-a|$ $\exists\,z_r\in\{z\in\mathbb{C}\,|\,|z-a|=r\}$
so that $f(z_r)=f(a)$. This will suffice, for we can take $a_n=z_{r_n}$ where $r_n=|b-a|/2^n$. So we fix $r$,
$0<r<|b-a|$. We recall that $\lim_{t\rightarrow\infty}f(g(\gamma(t)))=f(a)=f(b)$, $\lim g(\gamma(t_n))=a$,
$\lim g(\gamma(s_n))=b$, $\lim t_n=\lim s_n=\infty$. So the curve $g(\gamma(t))$, $t\in [0,\infty)$ goes back and forth
between points that are as close as we please to $a$ or to $b$, and for large $t$, $f$ maps this curve to a
neighborhood of $f(a)$ and this neighborhood can be as small as we please. Taking the intersection points of the 
curve $g(\gamma(t))$ with the circle $|z-a|=r$ (we recall that $0<r<|b-a|$) we see that $\exists$ a sequence
$\alpha_n$ on this circle such that $\alpha=\lim\alpha_n$ exists and $f(\alpha)=f(a)$. We take $z_r=\alpha$. $\qed $

\begin{remark}
We note that our proof of Proposition 3.6 uses the fact that $f\in {\rm elh}(\mathbb{C})$, but $g$ can be less restricted.
In fact $g$ can be only a continuous function.
\end{remark}
\noindent
An important consequence comes out of Proposition 3.6 and out of its proof. This will be a result that
sharpens Proposition 3.5. Namely,

\begin{theorem}
If $f,g\in {\rm elh}(\mathbb{C})$, then $A(f\circ g)=A(f)\cup f(A(g))$.
\end{theorem}
\noindent
{\bf Proof.} \\
By Proposition 3.5 we have $f(A(g))\subseteq A(f\circ g)$. If $a\in A(f\circ g)-f(A(g))$ then
there exists an asymptotic tract $\gamma\in T_f(f\circ g)-T_f(g)$ so that $\lim_{t\rightarrow\infty}f(g(\gamma(t)))=a$.
By Proposition 3.6 $\lim_{t\rightarrow\infty}g(\gamma(t))=\infty$. Hence $g(\gamma(t))\in T_f(f)$ and 
$a\in A(f)$. This proves that, $A(f\circ g)-f(A(g))\subseteq A(f)$. Finally, if $a\in A(f)-A(f\circ g)$ then
there is an asymptotic tract $\delta\in T_f(f)$, i.e. $\lim_{t\rightarrow\infty}\delta(t)=\infty$
and $\lim_{t\rightarrow\infty}f(\delta(t))=a$. But (since $a\not\in A(f\circ g)$) there is no 
$\gamma(t)\in T_f(f\circ g)$ for which $\lim_{t\rightarrow\infty}f(g(\gamma(t)))=a$. To get a contradiction
(and thus to prove that there is no $a\in A(f)-A(f\circ g)$) we consider the pre-image $g^{-1}(\delta)$.
Each component $\beta$ of this curve can not be bounded because $\delta\bigtriangleup g(\beta)$ is a finite
set and $\delta$ is unbounded. If $\lim_{t\rightarrow\infty}\beta(t)\ne\infty$ then as in the proof
of Proposition 3.6 we deduce that $f\circ g\equiv {\rm Const.}$. Hence not both $f,g\in {\rm elh}(\mathbb{C})$.
We conclude that $\lim_{t\rightarrow\infty}\beta(t)=\infty$ and $\lim_{t\rightarrow\infty}f(g(\beta(t)))=a$.
But this means that $a\in A(f\circ g)$ which contradicts our assumption that $a\in A(f)-A(f\circ g)$. $\qed $

\begin{remark}
We see that if $f,g\in {\rm elh}(\mathbb{C})$ and $a\in A(f)$ and $\delta\in T_f(f)$ is an $a$-asymptotic tract,
then for each irreducible component $\beta $ of the curve $g^{-1}(\delta)$ either $\beta\in T_f(g)$ or
$g(\beta(t))\in T_f(f)$ is an $a$-asymptotic tract of $f$.
\end{remark}

\begin{proposition}
Let $f\in {\rm elh}(\mathbb{C})$. If $\exists\,g\in {\rm elh}(\mathbb{C})$ such that $T_f(g)=T_f(f\circ g)$,
then $f(\mathbb{C})=\mathbb{C}$.
\end{proposition}
\noindent
{\bf Proof.} \\
Since $f\in {\rm elh}(\mathbb{C})$ we have $\mathbb{C}-f(\mathbb{C})\subseteq A(f)$ because the only
points in the complement of the image of $f$ are the Picard exceptional values of $f$ which are also
asymptotic values of $f$. Of course $f$ may have at most one such a value. If $T_f(g)=T_f(f\circ g)$ then 
by Theorem 3.8 we must have $A(f)\subseteq f(A(g))\subseteq f(\mathbb{C})$. Thus there are no
Picard exceptional values of $f$. $\qed $ \\
\\
The last proposition implies that if $f\in  {\rm elh}(\mathbb{C})$ is not a surjective mapping then 
$\forall\,g\in {\rm elh}(\mathbb{C})$ we must have $T_f(g)\subset T_f(f\circ g)$. In particular
$T_f(f)\subset T_f(f\circ f)$. But also $f\circ f$ is not surjective because $(f\circ f)(\mathbb{C})
\subseteq f(\mathbb{C})\subset\mathbb{C}$. Hence $T_f(f\circ f)\subset T_f(f\circ f\circ f)$. By induction
on the number of compositions we get the infinite strictly ascending sequence,
$$
T_f(f)\subset T_f(f\circ f)\subset T_f(f\circ f\circ f)\subset\ldots\subset T_f(f\circ f\circ\ldots\circ f)\subset\ldots.
$$
Another consequence of Theorem 3.8, is that for any $f\in {\rm elh}(\mathbb{C})$ we have the identity,
$$
A(f^{\circ(n+1)})=A(f)\cup\bigcup_{k=1}^n f^{\circ k}(A(f)).
$$
To see that we use induction as follows $A(f)\cup f(A(f))=A(f\circ f)$, hence $A(f)\cup f(A(f))\cup (f\circ f)(A(f))=
A(f)\cup f(A(f)\cup f(A(f)))=A(f)\cup f(A(f\circ f))=A(f\circ f\circ f)$, and so on.

\begin{proposition}
The following are equivalent: \\
{\rm 1)} $f\not\in {\rm Aut}(\mathbb{C})$. \\
{\rm 2)} $R_f({\rm elh}(\mathbb{C}))\subset {\rm elh}(\mathbb{C})-{\rm Aut}(\mathbb{C})$. \\
{\rm 3)} $L_f({\rm elh}(\mathbb{C}))\subset {\rm elh}(\mathbb{C})-{\rm Aut}(\mathbb{C})$.
\end{proposition}
\noindent
{\bf Proof.} \\
This follows by Theorem 3.8 and Proposition 3.5 which imply that:
$$
T_f(R_f(g))=T_f(g\circ f)\supseteq T_f(f),\,\,A(L_f(g))=A(f\circ g)\supseteq A(f),
$$
and by Proposition 2.2 part 5, which implies that $f\not\in {\rm Aut}(\mathbb{C})$ is equivalent
to $A(f)\ne\emptyset$ and so also to $T_f(f)\ne\emptyset$. $\qed $ 

\begin{proposition}
$\forall\,f\in {\rm elh}(\mathbb{C})$, $R_f$ is injective.
\end{proposition}
\noindent
{\bf Proof.} \\
$R_f(g)=R_f(h)\Rightarrow g\circ f=h\circ f\Rightarrow g|_{f(\mathbb{C})}=h|_{f(\mathbb{C})}\Rightarrow
g\equiv h$. The last implication follows by $|\mathbb{C}-f(\mathbb{C})|\le 1$. $\qed $

\begin{proposition}
$\forall\,f,g,h\in {\rm elh}(\mathbb{C})$, $(L_f(g)=L_f(h)\wedge g\not\equiv h)\Rightarrow\exists\,t(z)$
an entire function, such that $g(z)=h(z)+e^{t(z)}$.
\end{proposition}
\noindent
{\rm Proof.} \\
$L_f(g)=L_f(h)\Rightarrow f\circ g=f\circ h$. By $g\not\equiv h$ there are points $z\in\mathbb{C}$ for which
$g(z)\ne h(z)$. Let $N=\{z\in\mathbb{C}\,|\,g(z)\ne h(z)\}$. We will prove that $N=\mathbb{C}$ and so 
$g(z)-h(z)$ never vanishes, so $\exists$ an entire function $t(z)$ such that $g(z)-h(z)=e^{t(z)}$.
Thus $N\ne\emptyset$ is open in the strong topology. For $g,h$ are local homeomorphisms and if
$g(z)\ne h(z)$ then $\exists\,O$, an open neighborhood of $z$ in the strong topology, such that
$g(O)\cap h(O)=\emptyset$. Let $z\in\partial N$. Let $z_n\in N$ satisfy $\lim z_n=z$. Then 
$\forall\,n\in\mathbb{Z}^+,\,g(z_n)\ne h(z_n), f(g(z_n))=f(h(z_n))$ and $g(z)=h(z)$. This implies
that in any strong neighborhood of $g(z)=h(z)$ there are distinct points $g(z_n)\ne h(z_n)$ (for
$n\in\mathbb{Z}^+$ large enough), so that $f(g(z_n))=f(h(z_n))$. Hence $f$ can not be injective 
in any strong neighborhood of the point $g(z)=h(z)$. This contradicts the assumption $f\in {\rm elh}(\mathbb{C})$
and this proves that $\partial N=\emptyset$ or equivalently that $N=\mathbb{C}$. $\qed $ \\
\\
Surprisingly, ${\rm elh}(\mathbb{C})$ behaves very different from \'et$(\mathbb{C}^2)$, the two dimensional Keller mappings 
with respect to $L_f$. Namely in this case we have:

\begin{proposition}
$\exists\,f\in {\rm elh}(\mathbb{C})$, such $L_f$ is not injective.
\end{proposition}
\noindent
{\bf Proof.} \\
We take $f(z)=e^{2\pi iz}$, $g(z)=e^z+1$, $h(z)=e^z$. Then $g\not\equiv h$ but $L_f(g)=
e^{2\pi i(e^z+1)}=e^{2\pi ie^z}=L_f(h)$. $\qed $

\begin{remark}
Unless otherwise said our topology on ${\rm elh}(\mathbb{C})$ will be that of local uniform convergence.
Thus $f_n\rightarrow f$ in that topology if the sequence converges to $f$ uniformly on each compact subset
$K\subset\mathbb{C}$. We recall (by Proposition 2.2 part 2) that in this case $f\in {\rm elh}(\mathbb{C})$.
\end{remark}
\noindent
We recall the connection with the so-called compact-open topology. If $X$ and $Y$ are topological spaces, 
$C(X,Y)$ denotes the set of all the continuous mappings from $X$ to $Y$. If $K$ is a compact subset
of $X$ and $U$ an open subset of $Y$, we define $S(K,U)=\{f:\,X\rightarrow Y\in C(X,Y)\,|\,f(K)\subset U\}$.
The sets $S(K_1,\ldots,K_n;U_1,\ldots,U_n)=\bigcap_{i=1}^n S(K_i,U_i)$, $n\in\mathbb{Z}^+$, form a basis for a 
topology on $C(X,Y)$ called the compact-open topology. It is denoted by $\tau_{CO}$. The following facts are well 
known: \\
$1_{CO}$) In the compact-open topology, $C(X,Y)$ is a Hausdorff or regular space whenever $Y$ is a Hausdorff or
regular space. \\
$2_{CO}$) If $X$ is locally compact and $X$ and $Y$ are second countable, then so is $(C(X,Y),\tau_{CO})$. \\
$3_{CO}$) Let $X$ be locally compact and second countable, and let $Y$ be regular and second countable.
Then $(C(X,Y),\tau_{CO})$ is metrizable. \\
\\
Let $(Y,d)$ be a metric space and let $X$ be a topological space. For a compact subset $K\subset X$ of $X$,
$\varepsilon>0$, and $f\in Y^X$, we define $B_K(f,\varepsilon)=\{g:\,X\rightarrow Y\,|\,d(f(x),g(x))<\varepsilon,
\forall\,x\in K\}$. The sets $B_K(f,\varepsilon)$ form a basis for a topology on $Y^X$, called the topology
of compact convergence. It is denoted by $\tau_{CC}$. Clearly $f_n$ converges to $f$ in $\tau_{CC}$ if and only if
$\forall$ compact $K$, $f_n|_K$ converges uniformly to $f|_K$. \\
\\
$4_{CC-CO}$) $C(X,Y)$ is closed in $(Y^X,\tau_{CC})$. \\
$5_{CC-CO}$) Let $X$ be a topological space and $Y$ a metric space. Then $(C(X,Y),\tau_{CO})=(C(X,Y),\tau_{CC})$. \\
\\
In particular, by $5_{CC-CO}$ we have $({\rm elh}(\mathbb{C}),\tau_{CO})=({\rm elh}(\mathbb{C}),\tau_{CC})$. Our 
natural topology on ${\rm elh}(\mathbb{C})$ is $\tau_{CC}$, as mentioned before. Thus in this case it coincides
with the compact-open topology, $\tau_{CO}$.

\begin{proposition}
{\rm 1)} The topological space $({\rm elh}(\mathbb{C}),\tau_{CC})$ is path connected. \\
{\rm 2)} The image of the $f$-right mapping, $R_f({\rm elh}(\mathbb{C}))$ is a closed subset of
$({\rm elh}(\mathbb{C}),\tau_{CC})$. Here $f\in {\rm elh}(\mathbb{C})$ is fixed.
\end{proposition}
\noindent
{\bf Proof.} \\
1) Let us assume first that $f\in {\rm elh}(\mathbb{C})$ and also $f(0)=0$. Then $\forall\,0<t\le 1$,
$f_t(z)=(1/t)f(tz)\in {\rm elh}(\mathbb{C})$. Let us denoted $f_0(z)={\rm id}(z)$. Then $f_1=f$ and
$f_0={\rm id}$ and $\{ f_t\,|\,0\le t\le 1\}$ is a path in ${\rm elh}(\mathbb{C})$ from $f_0$ to $f$.
In the general case we take $(1/t)(f(tz)-f(0))+f(0)$ that connects $z+f(0)$ to $f(z)$. Now 
${\rm Aut}(\mathbb{C})$ is path connected and we are done. \\
2) By $3_{CO}$ the space $({\rm elh}(\mathbb{C}),\tau_{CO})$ is metrizable. By $5_{CC-CO}$ it
coincides with $({\rm elh}(\mathbb{C}),\tau_{CC})$ and so to prove the closedness of $R_f({\rm elh}(\mathbb{C}))$
we can show that the limit point of any convergent sequence in $R_f({\rm elh}(\mathbb{C}))$ belongs to
$R_f({\rm elh}(\mathbb{C}))$. That will be sufficient. Thus let $f_n\in R_f({\rm elh}(\mathbb{C}))$ satisfy
$f_n\rightarrow g$ in $\tau_{CC}$. By the definition of $R_f({\rm elh}(\mathbb{C}))$ $\forall\,n\in\mathbb{Z}^+,\,
\exists\,g_n\in {\rm elh}(\mathbb{C})$ such that $f_n=g_n\circ f$. we have $g_n\circ f\rightarrow g$ in $\tau_{CC}$.
So $\forall\,K\subset\mathbb{C}$ a compact and $\forall\,\varepsilon>0,\,\exists\,N=N(\varepsilon,K)$ such
that $n,m>N\Rightarrow |g_n(f(z))-g_m(f(z))|<\varepsilon\,\forall\,z\in K$. So $n,m>N\Rightarrow 
|g_n(w)-g_m(w)|<\varepsilon,\,\forall\,w\in K_1=f(K)$. This implies by Proposition 2.2 part 3 that
$h=\lim g_n\in {\rm elh}(\mathbb{C})$ and hence $f_n=g_n\circ f\rightarrow h\circ f\in R_f({\rm elh}(\mathbb{C}))$. $\qed $

\begin{remark}
The mapping $R_f:\,{\rm elh}(\mathbb{C})\rightarrow {\rm elh}(\mathbb{C})$ is continuous and an injective 
mapping (Proposition 3.12). The image $R_f({\rm elh}(\mathbb{C}))$ is a closed subset of $({\rm elh}(\mathbb{C}),\tau_{CC})$
(Proposition 3.16 part 2). Hence $R_f({\rm elh}(\mathbb{C}))$ is also open $\Leftrightarrow R_f({\rm elh}(\mathbb{C}))=
{\rm elh}(\mathbb{C})$ because by Proposition 3.16 part 1 the space $({\rm elh}(\mathbb{C}),\tau_{CC})$ is
connected. This is equivalent to $f\in {\rm Aut}(\mathbb{C})$ by Proposition 3.11.
\end{remark}
\noindent
Thus it is natural to ask the following: Let $f\in {\rm elh}(\mathbb{C})-{\rm Aut}(\mathbb{C})$. What is the
boundary of the image $\partial R_f({\rm elh}(\mathbb{C}))$? We note that since $R_f({\rm elh}(\mathbb{C}))$ is
closed in $({\rm elh}(\mathbb{C}),\tau_{CC})$, it follows that $\partial R_f({\rm elh}(\mathbb{C}))\subseteq R_f({\rm elh}(\mathbb{C}))$.
Moreover, every automorphism $z+a$ is an interior point of the open set ${\rm elh}(\mathbb{C})-R_f({\rm elh}(\mathbb{C}))$
because $z+a\not\in R_f({\rm elh}(\mathbb{C}))$. The reason is that if $z+a\in R_f({\rm elh}(\mathbb{C}))$ then
$g(f(z))=z+a$ for some $g\in {\rm elh}(\mathbb{C})$. This implies by Proposition 3.5 that $T_f(f)\subseteq T_f(z+a)$
which is an absurd since $T_f(f)\ne\emptyset$ while $T_f(z+a)=\emptyset$. So a boundary point 
$h\in\partial R_f({\rm elh}(\mathbb{C}))$ has the form $h=g\circ f$ for some $g\in {\rm elh}(\mathbb{C})$ but there 
should exist a sequence $h_n\in {\rm elh}(\mathbb{C})-R_f({\rm elh}(\mathbb{C}))$ such that $h_n\rightarrow h$
in $\tau_{CC}$. At most finitely many of the $h_n\in {\rm Aut}(\mathbb{C})$, because by the Theorem
of Hurwitz ${\rm Aut}(\mathbb{C})$ is closed in $({\rm elh}(\mathbb{C}),\tau_{CC})$. To try and obtain such
a sequence $h_n$, we might want to see if we can pick its elements from among the family $h_t(z)=(1/t)h(tz)$,
$0<t\le 1$. We assume here that $g(0)=f(0)=0$ for simplicity. We observe that $h_t(z)=(1/t)g(f(tz))$. 
So the question boils down to: can it possibly be that $\forall\,0<t<1$ we will have $h_t(z)\not\in R_f({\rm elh}(\mathbb{C}))$?
More concretely, can we find a family $G(t,f(z)),\,0<t<1$ of entire functions in $f(z)$ such that
$$
\frac{1}{t}g(f(tz))=G(t,f(z)),\,\,0<t\le 1.
$$
We observe that the left hand side is entire both in $t$ and in $z$, and in fact there is a full symmetry
between $t$ and $z$ in $g(f(tz))$ which indicates a symmetry in the power series coefficients (only the
diagonal terms can be different from zero). In other words $tG(t,f(z))=\sum_{n=1}^{\infty} A_{nn}t^n z^n$.
In fact $\forall\,n\in\mathbb{Z}^+$, $A_{nn}=(1/n!)(d^n/dz^n)\{g(f(z))\}|_{z=0}$. We can derive a first order
linear pde on $G(u,v)$ that avoids the entire function $g(z)$, as follows:
$$
g(f(tz))=tG(t,f(z)),
$$
applying $\partial_z$ we get,
$$
tf'(tz)g'(f(tz))=tf'(z)G_v(t,f(z)),
$$
hence,
$$
f'(tz)g'(f(tz))=f'(z)G_v(t,f(z)).
$$
On the other hand if we apply $\partial_t$ instead we get,
$$
zf'(tz)g'(f(tz))=G(t,f(z))+tG_u(t,f(z)).
$$
Taking the obvious linear combination of the two last equations we arrive at our pde,
$$
G(t,f(z))+tG_u(t,f(z))-zf'(z)G_V(t,f(z))=0.
$$
As promised the last equation does not contain $g$. We have noted before that $G(t,f(z))$ is
entire in $(tz)$ and so we deduce that $G(u,v)$ is entire in $u,v$. Let us assume that $G(u,v)=\sum_{n,m}B_{nm}u^nv^m$.
We substitute this into into our pde and obtain,
$$
0\equiv G(t,f)+tG_u(t,f)-zf'G_v(t,f)=
$$
$$
=\sum B_{nm}t^nf^m+\sum nB_{nm}t^nf^m-\frac{zf'}{f}\sum mB_{nm}t^nf^m=
$$
$$
=\sum\left(\left(1+n-\frac{zf'}{f}m\right)B_{nm}f^m\right)t^n.
$$
We deduce that $\forall\,n\in\mathbb{Z}^+\cup\{ 0\}$ we have the identities,
$$
\sum_m\left(1+n-\frac{zf'}{f}m\right)B_{nm}f^m\equiv 0.
$$
Let us fix an $n\in\mathbb{Z}^+\cup\{ 0\}$. Then
\begin{equation}
\label{eq1}
(1+n)\sum_m B_{nm}f^m\equiv (zf')\sum_m mB_{nm}f^{m-1}
\end{equation}
We denote $H_n(w)=\sum_m B_{nm}w^m$ and then $H'_n(w)=\sum_m mB_{nm}w^{m-1}$.
Then $(1+n)H_n(f)=(zf')H'_n(f)$ and hence,
$$
\left(\frac{1}{1+n}\right)\frac{H'_n(f)}{H_n(f)}=\frac{1}{zf'(z)}\,\,{\rm independent}\,{\rm of}\,n.
$$
So $\forall\,k,n\in\mathbb{Z}^+\cup\{ 0\}$, $(1/(1+n))(\log H_n(f))'=(1/(1+k))(\log H_k(f))'$ and we deduce
that $H^{1+k}_n=H^{1+n}_k$. Taking $k=0$ we get $H_n=H^{1+n}_0$. Substituting that into the identity (\ref{eq1})
we obtain $(1+n)H_0(f)^{1+n}=(zf')(1+n)H_0(f)^{1+n}$. So $H_0(f)=zf'$. Thus $H_n=H^{1+n}_0=(zf')^{1+n}$. We
recall that according to our notations we have $G(u,v)=\sum B_{nm}u^nv^m$, and hence
$$
G(t,f)=\sum_{n,m}B_{nm}u^nv^m|_{(u,v)=(t,f)}=\sum_{n,m} B_{nm}t^nf^m=\sum_n(\sum_m B_{nm}f^m)t^n=
$$
$$
=\sum_n H_n(f)t^n=\sum_n(zf')^{1+n}t^n=(zf')\sum_n(tzf'(z))^n=\frac{zf'}{1-tzf'}.
$$
We strongly use here the permanence principle. We deduce that
$$
\frac{(tz)f'(z)}{1-(tz)f'(z)}=g(f(tz)).
$$
In the last step we used $G(t,f)=(1/t)g(f(tz))$. In particular for $t=1$ we obtain the following 
functional equation,
$$
\frac{zf'(z)}{1-zf'(z)}=g(f(z)).
$$
Since the right hand side is an entire function we deduce that $\forall\,z\in\mathbb{C}$, $zf'(z)\ne 1$. We
 now recall a result from Nevanlinna's value distribution theory for meromorphic functions. \\
 \\
 {\bf Theorem (R. Nevanlinna).} {\it Let $q(z)$ be a transcendental meromorphic function. Then for
 each value $a$, finite or infinite, the equation $q(z)=a$ has an infinite number of roots except
 for at most two exceptional values.} \\
 \\
 Using that with $q(z)=zf'(z)$ we see that for each of $a=0,\infty $ the equation has exactly one root.
 For $a=1$ it has no roots and we deduce that $zf'(z)$ is not a transcendental meromorphic function.
 This contradiction shows that our original equation 
 $$
 \frac{1}{t}g(f(tz))=G(t,f(z)),\,\,0<t\le 1,
 $$
 is contradictory and so we can find a sequence of real numbers $t_n\rightarrow 1^-$ such that
 $\forall\,n\in\mathbb{Z}^+$,
 $$
 h_{t_n}(z)=\frac{1}{t_n}g(f(t_nz))\not\in R_f({\rm elh}(\mathbb{C})).
 $$
 We have $\lim_{n\rightarrow\infty} h_{t_n}(z)=h(z)=g(f(z))\in R_f({\rm elh}(\mathbb{C}))$. This proves that
 $R_f({\rm elh}(\mathbb{C}))\subseteq\partial R_f({\rm elh}(\mathbb{C}))$ and so we have proved the following 
 interesting,
 \begin{theorem}
 $\forall\,f\in {\rm elh}(\mathbb{C})-{\rm Aut}(\mathbb{C}),\,\,\,\partial R_f({\rm elh}(\mathbb{C}))=
 R_f({\rm elh}(\mathbb{C}))$.
 \end{theorem}
 \noindent
 It is natural to ask if for all $f\in {\rm elh}(\mathbb{C})$, the image $L_f({\rm elh}(\mathbb{C}))$ is a closed 
 subset of $({\rm elh}(\mathbb{C}),\tau_{CC})$. The key for proving the parallel claim for $R_f({\rm elh}(\mathbb{C}))$
 (Proposition 3.16 part 2) was the tame behavior of $\tau_{CC}$-convergent sequences. Namely, if
 $f,g_n\in {\rm elh}(\mathbb{C})$ and $\{g_n\circ f\}$ is $\tau_{CC}$-convergent then $\{ g_n\}$ is
 $\tau_{CC}$-convergent. This tameness is false for $L_f$. For example, $f(z)=e^z$, $g_n(z)=z+2\pi in$ satisfy
 by the periodicity of $f$, $f\circ g_n=f$ but $\{ g_n\}$ is not $\tau_{CC}$-convergent. Proposition 3.14
 tells us that $L_f$ might not be injective. With the aid of Proposition 3.13 we can make this much
 more precise.
 
 \begin{proposition}
 $\forall\,f,g\in {\rm elh}(\mathbb{C})$ we have $L^{-1}_f(L_f(g))=\{g(z)+k_je^{t(z)}\,|\,j=0,\ldots,N,\,k_0=0\}$,
 where $N\in\mathbb{Z}^+\cup\{ 0,\infty\}$. Moreover, we have,
 $$
 N\,{\rm is}\,{\rm finite}\,\Leftrightarrow \left\{\begin{array}{ll} L^{-1}_f(L_f(g))=\{g\} & (N=0,\,k_0=0), \\
 {\rm or} &  \\ |L^{-1}_f(L_f(g))|>1\,{\rm in}\,{\rm which}\,{\rm case}: & f\in\mathbb{C}[z],\,N=\deg f-1.\end{array}\right.
 $$
 \end{proposition}
 \noindent
 {\bf Proof.} \\
 Clearly $g\in L^{-1}_f(L_f(g))$ and by Proposition 3.13 if $h\ne g$, $h\in L^{-1}_f(L_f(g))$ then 
 $h(z)=g(z)+e^{t(z)}$ for some entire $t(z)$. Let us assume further that $h_1(z)\in L^{-1}_f(L_f(g))-\{ g,h\}$. Then
 by Proposition 3.13 there are two entire functions $t_1(z),\,t_2(z)$ such that $h_1(z)=g(z)+e^{t_1(z)}$,
 $h_1(z)=h(z)+e^{t_2(z)}$. Hence $g(z)+e^{t_1(z)}=(g(z)+e^{t(z)})+e^{t_2(z)}$, and so we obtain
 $1\equiv e^{t(z)-t_1(z)}+e^{t_2(z)-t_1(z)}$. This implies that the entire functions $e^{t(z)-t_1(z)}$ and
 $e^{t_2(z)-t_1(z)}$ do not assume the values $\{ 0,1\}$ and hence are constant by Picard's Theorem.
 Thus $t_1(z)=t(z)+c_1$, $t_2(z)=t(z)+c_2$ for some $c_1,c_2\in\mathbb{C}$ which proves that any function
 in $L^{-1}_f(L_f(g))$ is of the form $g(z)+ke^{t(z)}$, $k\in\mathbb{C}$. Since $\mathbb{C}$ is tiled up
 by a countable set of pairwise disjoint maximal domains of $f$ the set of admissible constants $k$ is a 
 countable set, say, $\{ k_j\,|\,j=0,\ldots,N,\,k_0=0\}$. Note that if $N\ge 1$ we may assume that $k_1=1$.
 The only entire functions with finitely many pairwise disjoint maximal domains are polynomials, \cite{rp}, $f\in\mathbb{C}[z]$
 and the number of tiles for such an $f(z)$ is $\deg f$. Thus in this case $N=\deg f -1$. $\qed $
 
 \begin{proposition}
 The set $\{ k_je^{t(z)}\,|\,j\in I,\,k_0=0\}$ is a cyclic subgroup of $(\mathbb{C},+)$.
 \end{proposition}
 \noindent
 {\bf Proof.} \\
 Given $f,g\in {\rm elh}(\mathbb{C})$ we saw in Proposition 3.19 that
 $$
 L^{-1}_f(L_f(g))=\{g(z)+k_je^{t(z)}\,|\,j\in I,\,k_0=0\},\,\,|I|\le\aleph_0.
 $$
 Thus for a generic $z$, the points $g(z)+k_j e^{t(z)}$ are the full set of points that are $f$-equivalent to $g(z)$.
 However, for each $g(z)$ the set of all $f$-equivalent points to $g(z)$ are the orbit of $g(z)$ under
 the action of the group of the $f$-deck transformations. It is clear that the composition of the two
 $f$-deck transformations $g(z)\rightarrow g(z)+k_1 e^{t(z)}$, $g(z)\rightarrow g(z)+k_2 e^{t(z)}$ is
 $g(z)\rightarrow g(z)+(k_1+k_2)e^{t(z)}$. $\qed $
 
 \begin{proposition}
 If $f\in\mathbb{C}[z]$, then $L_f$ is injective. 
 \end{proposition}
 \noindent
 {\bf Proof.} \\
 By Proposition 3.19 and Proposition 3.20 in this case $\{k_je^{t(z)}\,|\,j\in I,\,k_0=0\}$ is a finite cyclic
 subgroup of $(\mathbb{C},+)$ and hence trivial. $\qed $
 
 \begin{remark}
 The last statement in the proof of Proposition 3.20 follows by: Let $h_(z),h_2(z)\in L^{-1}_f(L_f(g))$. Then we note
 that we have the two identities $L^{-1}_f(L_f(g))=L^{-1}_f(L_f(h_1))=L^{-1}_f(L_f(h_2))$. As explained in the proof of
 Proposition 3.19 there is an entire $t(z)$ and two constants $k_1,k_2$ so that $h_1(z)=g(z)+k_1 e^{t(z)}$,
 $h_2(z)=g(z)+k_2 e^{t(z)}$. Also by $h_2\in L^{-1}_f(L_f(h_1))$ there exists a constant $k_3$ so that
 $h_2(z)=h_1(z)+k_3 e^{t(z)}$ (same entire function $t(z)$). We can use the symbols $k_1,k_2$ and $k_3$ to denote
 the corresponding deck-transformations of $f^{-1}\circ f$. Indeed by $h_2(z)=h_1(z)+k_3 e^{t(z)}$ we have
 $g(z)+k_2 e^{t(z)}=(g(z)+k_1 e^{t(z)})+k_3e^{t(z)}=g(z)+(k_1+k_2)e^{t(z)}$, so the composition
 $k_2=k_3\circ k_1$ of the deck-transformations corresponds to complex addition $k_2=k_3+k_1$.
 \end{remark}
 
 \begin{remark}
 If $f\in {\rm elh}(\mathbb{C})-\mathbb{C}[z]$, $g\in {\rm elh}(\mathbb{C})$ and $|L^{-1}_f(L_f(g))|>1$, then
 by Proposition 3.19 we have $|L^{-1}_f(L_f(g))|=\aleph_0$ so by Proposition 3.20 the set $\{ k_j e^{t(z)}\,|\,j\in I,\,k_0=0\}$ 
 is an infinite cyclic group for addition. But then by changing, if necessary, the entire $t(z)$ we may
 assume that the set is, in fact $\{je^{t(z)}\,|\,j\in\mathbb{Z}\}$. 
 \end{remark}
 
 \begin{proposition}
 $\forall\,f\in {\rm elh}(\mathbb{C})-\mathbb{C}[z]$, $g\in {\rm elh}(\mathbb{C})$ if $|L^{-1}_f(L_f(g))|>1$ then
 there exists an entire $t(z)$ in the variable $z$ and an entire $h(z,w)$ in $(z,w)$ so that
 $$
 f(g(z)+we^{t(z)})=f(g(z))+e^{h(z,w)}\sin\pi w.
 $$
 \end{proposition}
 \noindent
 {\bf Proof.} \\
 We consider $f(g(z)+we^{t(z)})-f(g(z))$ with the entire function $t(z)$ of Proposition 3.19 normalized
 as in the last remark. Then this function is entire in $(z,w)$ and for a fixed $z$ it vanishes exactly
 for $w=j\in\mathbb{Z}$. Since$\sum_{j=1}^{\infty}(1/j)=\infty$, $\sum_{j=1}^{\infty}(1/j^2)<\infty$ we
 get by writing the Weierstrass canonical infinite product for our function,
 $$
 f(g(z)+we^{t(z)})-f(g(z))=e^{h_0(z,w)}w\prod_{j\ne 0,j=-\infty}^{\infty}\left(1-\frac{w}{j}\right)e^{w/j}=
 $$
 $$
 =e^{h_0(z,w)}w\prod_{j=1}^{\infty}\left(1-\frac{w^2}{j^2}\right)=\frac{1}{\pi}e^{h_0(z,w)}\sin\pi w,
 $$
 where in the last step we used the infinite product representation of $\sin\pi w$ which extends the 
 Euler-Wallis' formula for $\pi$:
 $$
 \frac{\pi}{2}=\frac{2}{1}\cdot\frac{2}{3}\cdot\frac{4}{3}\cdot\frac{4}{5}\cdot\frac{6}{5}\cdot\frac{6}{7}\cdot\ldots=
 \prod_{j=1}^{\infty}\left(\frac{4j^2}{4j^2-1}\right).
 $$
 $\qed $
 
 \begin{example}
 $f(z)=e^{2\pi iz},\,\,g(z)=e^z,\,\,L^{-1}_f(L_f(g))=\{e^z+j\,|\,j\in\mathbb{Z}\}=\{e^z+je^0\,|\,j\in\mathbb{Z}\}$.
 So in this case the normalized entire $t(z)\equiv 0$, and $f(g(z)+we^{t(z)})=f(e^z+w)=e^{2\pi i(e^z+w)}=
 e^{2\pi ie^z}e^{2\pi iw}=e^{2\pi ie^z}(1+(e^{2\pi i w}-1))=f(g(z))+e^{2\pi ie^z}(e^{2\pi iw}-1)=f(g(z))+
 e^{2\pi ie^z}[(\cos 2\pi w-1)+i(\sin 2\pi w]=f(g(z))+e^{2\pi ie^z}[-2\sin^2\pi w+i2\cos\pi w\sin\pi w]=
 f(g(z))+2ie^{2\pi ie^z}(\cos\pi w+i\sin\pi w)\sin\pi w=f(g(z))+2ie^{\pi i(2e^z+w)}\sin\pi w=
 f(g(z))+2ie^{\pi i(2e^z+w)}\sin\pi w=f(g(z))+e^{\log 2+\pi i((1/2)+2e^z+w)}\sin\pi w$. Thus in this case
 $t(z)\equiv 0$ and $h(z,w)=\log 2+\pi i((1/2)+2e^z+w)$.
 \end{example}
 
 \begin{proposition}
 Let $f\in {\rm elh}(\mathbb{C})-\mathbb{C}[z]$, $g\in {\rm elh}(\mathbb{C})$ satisfy $|L^{-1}_f(L_f(g))|>1$ and
 let $t(z)$, $h(z,w)$ be the entire functions for which
 $$
 f(g(z)+we^{t(z)})=f(g(z))+e^{h(z,w)}\sin\pi w\,\,\,\,\,(\rm Proposition)\,\,3.24.
 $$
 Then there exists an entire $k(z,w)$ in $(z,w)$ so that
 $$
 \frac{\partial h}{\partial w}(z,w)\cdot\sin\pi w+\pi\cos\pi w=e^{k(z,w)}.
 $$
 \end{proposition}
 \noindent
 {\bf Proof.} \\
 We differentiate with respect to $w$ the functional equation $f(g(z)+we^{t(z)})=f(g(z))+e^{h(z,w)}\sin\pi w$.
 We obtain
 $$
 e^{t(z)}f'(g(z)+we^{t(z)})=\frac{\partial h}{\partial w}(z,w)e^{h(z,w)}\sin\pi w+\pi e^{h(z,w)}\cos\pi w,
 $$
 $$
 \frac{\partial h}{\partial w}(z,w)\sin\pi w+\pi\cos\pi w=e^{t(z)-h(z,w)}f'(g(z)+we^{t(z)}).
 $$
 Since $f\in {\rm elh}(\mathbb{C})$ it follows that $f'$ never vanishes. Hence the entire function
 $$
 \frac{\partial h}{\partial w}(z,w)\cdot\sin\pi w+\pi\cos\pi w,
 $$
 never vanishes and we are done. $\qed $
 
 \begin{remark}
 Proposition 3.21 has a truly trivial proof. For by the Fundamental Theorem of the Algebra $f\in\mathbb{C}[z]\cap 
 {\rm elh}(\mathbb{C})\Leftrightarrow f\in{\rm Aut}(\mathbb{C})$.
 \end{remark}
 \noindent
 we can sharpen the statement of Proposition 3.24 as follows,
 
 \begin{proposition}
 $\forall\,f\in {\rm elh}(\mathbb{C})-\mathbb{C}[z]$, $g\in {\rm elh}(\mathbb{C})$ if $|L^{-1}_f(L_f(g))|>1$ then
 there exist three entire functions $t(z),\,h(z)$ and $L(w)$ so that
 $$
 f(g(z)+we^{t(z)})=f(g(z))+e^{L(w)+h(z)}\sin\pi w.
 $$
 \end{proposition}
 \noindent
 {\bf Proof.} \\
 we need to prove that the entire function $h(z,w)$ has the structure $h(z,w)=L(w)+h(z)$ (we abuse the notation $h$).
 By Proposition 3.26 we have
 $$
 \frac{\partial h}{\partial w}(z,w)\cdot\sin\pi w+\pi\cos\pi w=e^{k(z,w)}.
 $$
 Let $L(z,w)=(1/\pi)(\partial h/\partial w)$. Then we conclude that $L(z,w)\sin\pi w+\cos\pi w$ 
 never vanishes. By changing our notation for $k(z,w)$ we conclude that $L(z,w)\sin\pi w+\cos\pi w=
 e^{k(z,w)}\Rightarrow L_z\sin\pi w=k_ze^k$. Let us assume in order to get a contradiction
 that $L_z\not\equiv 0$. Then also $k_z\not\equiv 0$ and so $(L_z/k_z)\sin\pi w=e^k$. Thus
 $$
 \frac{L_z}{k_z}\sin\pi w=L\sin\pi w+\cos\pi w,
 $$
 $$
 \Rightarrow\left\{\begin{array}{rrrrr} (L-L_z/k_z)\sin\pi w & + & \cos\pi w & \equiv & 0 \\
 \sin\pi w & + & \cos\pi w & = & e^{\theta(w)} \end{array}\right.
 $$
 The second equation is a consequence of $\sin^2\pi w+\cos^2\pi w\equiv 1$. We can consider the
 last system as a linear system in the unknowns $\sin\pi w$ and $\cos\pi w$. This system is non-homogeneous,
 the equations are independent and it is consistent $\forall\,(z,w)\in\mathbb{C}^2$. Hence the
 coefficients matrix never vanishes. Thus
 $$
 \left|\begin{array}{cc} L-L_z/k_z & 1 \\ 1 & 1 \end{array}\right|=e^{k_1(z,w)},
 $$
 $$
 \Rightarrow \left(L-\frac{L_z}{k_z}\right)-1=e^{k_1}.
 $$
 Multiplying the last equation by $\sin\pi w$ gives
 $$
 \left(L-\frac{L_z}{k_z}\right)\sin\pi w-\sin\pi w=e^{k_1}\sin\pi w,
 $$
 $$
 \Rightarrow L\sin\pi w-\sin\pi w-\frac{L_z}{k_z}\sin\pi w=e^{k_1}\sin\pi w,
 $$
 $$
 \Rightarrow L\sin\pi w-\sin\pi w-(L\sin\pi w+\cos\pi w)=e^{k_1}\sin\pi w,
 $$
 $$
 \Rightarrow -\sin\pi w-\cos\pi w=e^{k_1}\sin\pi w.
 $$
 This is a contradiction as the substitution $w=0$ shows. This shows that our assumption $L_z\not\equiv 0$
 is contradictory. Hence $L_z(z,w)\equiv 0$ and so $L(z,w)=L(w)$ depends on $w$ only. Thus 
 $$
 \frac{\partial h}{\partial w}(z,w)=\pi L(w).
 $$
 So $h(z,w)=L(w)+h(z)$. $\qed $
 
 \begin{remark}
 The equation in the last proposition can be written as follows,
 $$
 e^{-h(z)}f(g(z)+we^{t(z)})=e^{-h(z)}f(g(z))+e^{L(w)}\sin\pi w.
 $$
 \end{remark}
 \noindent
 If we differentiate this with respect to $w$ we get,
 $$
 e^{t(z)-h(z)}f'(g(z)+we^{t(z)})=(L'(w)\sin\pi w+\pi\cos\pi w)e^{L(w)}.
 $$
 Since $f'$ never vanishes (because $f\in {\rm elh}(\mathbb{C})$) we deduce (in accordance to the
 proof of Proposition 3.28) that
 $$
 L'(w)\sin\pi w+\pi\cos\pi w=e^{k(w)}.
 $$
 In fact $e^{k(w)}=c\cdot e^{-L(w)}$ because the two sides of the equation above must be constant (the
 left hand side depends only on $z$ while the right hand side depends only on $w$). If we differentiate 
 the identity in Remark 3.29 with respect to $z$ we get,
 $$
 -h'f(g+we^t)+(g'+wt'e^t)f'(g+we^t)=-h'f(g)+g'f'(g).
 $$
 This proves the following,
 
 \begin{proposition}
 Under the condition of Proposition 3.28 we have: \\
 {\rm a)} $L'(w)\sin\pi w+\pi\cos\pi w$ never vanishes. \\
 {\rm b)} The function $-h'(z)f(g(z)+we^{t(z)})+(g'(z)+wt'(z)e^{t(z)})f'(g(z)+we^{t(z)})$ is independent of $w$.
 \end{proposition}
 \noindent
 If in part (b) we substitute $w=-g(z)e^{-t(z)}$ and recall that $f(0)=0$, $f'(0)=1$ we get after the
simplification $g'(z)-g(z)t'(z)$. We conclude the following,

\begin{proposition}
$-h'(z)f(g(z))+g'(z)f'(g(z))\equiv g'(z)-g(z)t'(z)$.
\end{proposition}
\noindent
We can do better by differentiating the function in Proposition 3.30 part (b)with respect to $w$,

\begin{proposition}
$(t'(z)-h'(z))f'(g(z)+we^{t(z)})+(g'(z)+wt'(z)e^{t(z)})f''(g(z)+we^{t(z)})\equiv 0$.
\end{proposition}

\begin{theorem}
Let $f\in {\rm elh}(\mathbb{C})$. Then $L_f$ is not injective if and only if
$$
f(z)=\frac{1}{b}e^{bz}+a\,\,\,{\rm for}\,\,{\rm some}\,\,\,a\in\mathbb{C},\,b\in\mathbb{C}^{\times}.
$$
\end{theorem}
\noindent
{\bf Proof.} \\
Clearly $L_{e^{bz}/b+a}$ is not injective. In the other direction, if $L_f$ is not injective
we have the identity of Proposition 3.32. First, we claim that $t'(z)\equiv 0$. If not, then we can
substitute $w=-g'(z)e^{-t(z)}/t'(z)$ into that identity. Since $f'$ never vanishes we conclude that
$t'(z)-h'(z)\equiv 0$. Hence $(g'(z)+wt'(z)e^{t(z)})f''*g(z)+we^{t(z)})\equiv 0$. Now $f''(g(z)+we^{t(z)})\not\equiv 0$
because $w$ and $z$ are independent and $e^{t(z)}\ne 0$ and so in that case we had to conclude
that $f''(z)\equiv 0$. But then $f\in\mathbb{C}[z]\cap {\rm elh}(\mathbb{C})$ so by Proposition 3.21
$L_f$ is injective. Thus $g'(z)+wt'(z)e^{t(z)}\equiv 0$ which in turn (if $w=0$) implies that
$g'(z)\equiv t'(z)\equiv 0$. This contradicts our assumption that $t'(z)\not\equiv 0$. This proves that 
$t'(z)\equiv 0$ so that $t(z)\equiv c$ a constant. We get from the identity in Proposition 3.32,
$$
-h'(z)f'(g(z)+we^c)+g'(z)f''(g(z)we^c)\equiv 0.
$$
We now substitute $w=-g(z)e^{-c}$ and obtain $h'(z)f'(0)+g'(z)f''(0)\equiv 0$. We recall that
$f'(0)=1$ and get $g(z)=\alpha h(z)+\beta$, $\alpha,\beta\in\mathbb{C}$. So
$$
-h'(z)f'(g(z)+we^c)\alpha h'(z)f''(g(z)+we^c)\equiv 0.
$$
Now $g(z)\in {\rm elh}(\mathbb{C})$ and so $g'(z)\not\equiv 0$ and $h'(z)\not\equiv 0$. Thus,
$$
-f'(g(z)+we^c)+\alpha f''(g(z)+we^c)\equiv 0.
$$
We conclude that $f'(z)=\alpha f''(z)$ for some $\alpha\in\mathbb{C}^{\times}$. Hence $f(z)=\alpha f'(z)+\delta$
and so $f(z)=e^{bz}/b+a$. $\qed $

\begin{remark}
Theorem 3.33 tells us that the only obstacle for $L_f$ to be injective for any $f\in {\rm elh}(\mathbb{C})$ is that
$f$ is a translation of a conjugation by a $b\in\mathbb{C}^{\times}$ of the exponential function, $e^z$.
\end{remark}

\begin{remark}
Let us consider again Proposition 3.30 part (b), with $t(z)\equiv c$ a constant. Then the following is 
independent of $w$:
$$
-h')z)f(g(z)+we^c)+g'(z)f'(g(z)+we^c).
$$
At this point we know that $f(z)=e^z+a$ (we took $b=1$ for simplicity). So
$$
-h'(z)(e^{g(z)+we^c}+a)+g'(z)f'(g(z)+we^c),
$$
is independent of $w$. The above equals
$$
e^{we^c}(-h'(z)e^{g(z)}+g'(z)e^{g(z)})-h'(z)a,
$$
and we conclude that $-h'(z)+g'(z)\equiv 0$ which is consistent with the proof of Theorem 3.33.
\end{remark}

\begin{remark}
For later needs we point to the fact that, unfortunately, the family ${\rm elh}(\mathbb{C})$ is not a normal
family. For example we may consider the orbit of $e^z$ under conjugation with $\mathbb{Z}^+$, i.e. the sequence
$$
\left\{\frac{1}{n}e^{nz}\right\}^{\infty}_{n=1}.
$$
This sequence is contained in ${\rm elh}(\mathbb{C})$ but has no convergent subsequence (uniformly on compacta).
For if $z=x\in\mathbb{R}^+$ then $\lim_{n\rightarrow\infty}e^{nz}/n=\infty$ while if $z=x\in\mathbb{R}^-\cup\{0\}$
then $\lim_{n\rightarrow\infty}e^{nz}/n=0$ and for $y\in\mathbb{R}$ which is rationally independent of $\pi $ the
numbers $e^{\log n+iy}/n$ form a dense subset of $|z|=1$ (in the standard topology).
\end{remark}

\noindent
{\it Ronen Peretz \\
Department of Mathematics \\ Ben Gurion University of the Negev \\
Beer-Sheva , 84105 \\ Israel \\ E-mail: ronenp@math.bgu.ac.il} \\ 
 
\end{document}